\documentclass[11pt,a4paper]{amsart}

\usepackage[german,english]{babel}
\usepackage{amssymb}
\usepackage{amsmath}
\usepackage{pdfpages}
\usepackage{mathrsfs}
\usepackage{stmaryrd} 

\newcommand{\dd}{\mathrm{d}}
\newcommand{\II}{\mathrm{I\!I}}
\newcommand{\llangle}{\left\langle}
\newcommand{\rrangle}{\right\rangle}

\newcommand{\deltatnought}{\left.\frac{\partial}{\partial t}\right\vert_{t=0}}
\newcommand{\Order}{\mathscr{O}}
\newcommand{\grad}{\mathrm{grad}} 
\newcommand{\tangentpart}{\mathfrak{t}}

\newtheorem{theorem_}{Theorem}
\newtheorem{theoremofbeltrami_}[theorem_]{Theorem of Beltrami}
\newtheorem{lemma_}[theorem_]{Lemma}

\newtheorem{definition_}[theorem_]{Definition}

\begin{document}

\title{Some Modifications of the Theorem of Beltrami.}

\author{Steven Verpoort}

\thanks{
      \begin{tabular}{ll}
              Steven Verpoort. & \\
              \emph{Former address:} &
              \emph{Present address:} \\
              Departement Wiskunde &
              \'Ustav Matematiky a Statistiky \\
              Katholieke Universiteit Leuven &
              Masarykova Univerzita\\
              Celestijnenlaan 200B &
              Kotl\'a\v{r}sk\'a 2\\
              3001 Heverlee &
              611 37 Brno\\
              Belgium.  &
              Czech Republic.\\
       \end{tabular}      
              }

\maketitle

\begin{abstract}
The two main topics of this text are as follows:
Firstly, three modifications of the theorem of Beltrami will be presented for diffeomorphisms between Riemannian manifolds and a space form which preserve the geodesic circles, the geodesic hyperspheres, or the minimal surfaces, respectively.
Secondly, it is defined what it means for an infinitesimal deformation of a metric to preserve the geodesics up to first order, and a corresponding infinitesimal version of Beltrami's theorem is given.
\end{abstract}

\setcounter{tocdepth}{2}
\begin{center}
\rule{0.6\textwidth}{1.0pt}
\tableofcontents
\rule{0.6\textwidth}{1.0pt}
\end{center}

\newpage

\section*{Introduction.}

In this article will be presented
some modifications of the following classical theorem of E.\ Beltrami 
according to which the space forms are distinguished among all other Riemannian manifolds by the fundamental property of admitting local diffeomorphisms to a Euclidean space under which the geodesics correspond \cite{beltrami1865,beltrami1868}.

\begin{theoremofbeltrami_}
A Riemannian manifold which can be mapped onto a space of constant Riemannian curvature with correspondence of geodesics, has itself constant Riemannian curvature.

Consequently each of these two manifolds can be locally embedded as
either
a hyperplane not passing through the origin in a Euclidean space,
or
a hypersphere centred at the origin in a Euclidean space,
or
the standard embedding of a hyperbolic space in a Lorentz-Minkowski space as a hyperquadric centred at the origin. 
Now a diffeomorphism between these manifolds preserves the geodesics if and only if it originates by central projection from an affine transformation between these pseudo-Euclidean spaces which maps the centre to the centre.
\end{theoremofbeltrami_}

The question, whether the above theorem remains valid for diffeomorphism between Riemannian manifolds which preserve some other naturally defined classes of submanifolds, is a first main topic of this text. Therefore it seems of interest to introduce the following uniformous terminology: a diffeomorphism between two Riemannian manifolds will be said to be
\begin{itemize}
\item[(i).]
\textit{cogeodesical}, if it preserves the geodesics. This name should be compared with the nomenclature \textit{collineations} from classical projective geometry. We will not use the terms ``geodesic'' or ``projective;''
\item[(ii).]
\textit{concircular}, if it preserves the geodesic circles;
\item[(iii).]
\textit{cospherical}, if it preserves the geodesic hyperspheres;
\item[(iv).]
\textit{cominimal}, if it preserves the minimal hypersurfaces.
\end{itemize}

\textit{In the first part} of this article, some known facts on concircular diffeomorphisms will be presented as the concircular theorem of Beltrami.

\textit{In the second part,} it will be shown that a Riemannian manifold which  admits a cospherical diffeomorphism onto a space of constant Riemannian curvature, has itself constant Riemannian curvature. Moreover, for spaces of constant Riemannian curvature, the class of cospherical diffeomorphisms coincides with the class of concircular diffeomorphisms.

\textit{In the third part,} cominimal diffeomorphisms between three-dimensional Riemannian manifolds are studied. It should be remarked that the class of minimal surfaces in a given three-dimensional Riemannian manifold
depends on an \textit{infinity} of parameters, whereas the three other classes merely constitute a finite-dimensional class of submanifolds.
Therefore, prescribing all minimal surfaces for a Riemannian metric on a three-dimensional manifold seems a very severe restriction. It will be shown indeed that except for homotheties there do not exist cominimal diffeomorphisms.

\textit{In the last part} will be given an infinitesimal version of Beltrami's theorem, which is the second main topic of the text. First, the concept of 
an \textit{infinitesimally cogeodesical deformation} is defined, and a lemma which analytically translates this geometrical definition is presented 
(\S\,\ref{sec:def_infproj},\,\ref{sec:lem_infproj}). From this we can already prove the intended infinitesimal theorem of Beltrami in the  two-dimensional case (\S\,\ref{sec:two_dim}); in fact, the proof which is presented here is merely a small adaption of the proof of the classical theorem as it appears in, \textit{e.g.}, \cite{struik}. As such this proof is merely a technical calculation, which is straightforward but not well-suited for a higher-dimensional generalisation. 

A proof of the same result which holds for any dimension is essentially based on the known fact that the equation which describes all metrics which share their geodesics with a given metric, can be rewritten as a linear equation. This bridge between the ``finite'' and the ``infinitesimal'' problem of cogeodesical equivalence is the key element of the proof of the corresponding theorem for arbitrary dimensions (thm.~\ref{thm:inf_beltrami_arbdim} in \S\,\ref{sec:arb_dim}).

\section*{Notation and Assumptions.}

All manifolds are assumed to be connected.
We will say that a Riemannian manifold has \textit{constant Riemannian curvature} if the sectional curvature along a tangent plane is a constant which depends neither on the footpoint nor on the direction of the chosen tangent plane. We will use the term \textit{pointwise constant sectional curvature} if, at every fixed point, it is independent of the direction.

We will use the following sign convention concerning the Riemann curvature tensor of a Riemannian manifold $(M,g)$: if $\nabla$ stands for the Levi-Civita connection, then
$R(V,W)X=\nabla_{[V,W]}X-\nabla_V \left(\nabla_W X\right) + \nabla_W\left( \nabla_V X\right)$.

The \textit{Hessian operator} of a function $\varphi:M\rightarrow \mathbb{R}$ on a Riemannian manifold $(M,g)$ is defined as
\[
\textrm{Hs}_{\varphi}: \mathfrak{X}(M) \rightarrow \mathfrak{X}(M) : V \mapsto \nabla_V (\grad\varphi)\,.
\]
Here we have denoted $\mathfrak{X}(M)$ for the set of all vector fields on $M$. The directional derivative of a function $f$ along a vector field $V$ will be denoted by
$V[f]$.

\section{Diffeomorphisms which Preserve Geodesic Circles.}

By definition, a \textit{geodesic circle} in a Riemannian manifold is a curve for which the first geodesic curvature is constant and the second geodesic curvature vanishes. Diffeomorphisms which preserve geodesic circles, or \textit{concircular diffeomorphisms}, have been studied by, \textit{a.o.}, A. Fialkow, W. Vogel and K. Yano. The theorem below follows from \cite{kuhnel1988}, \S\,C.

\begin{theorem_}[``\textit{Concircular Beltrami theorem}'']
If a concircular diffeomorphism $\Psi$ between two Riemannian manifolds, one of which has constant Riemannian curvature, exists, then this diffeomorphism is conformal and both spaces have constant Riemannian curvature.

Assume now conversely that $\Psi:(M,g)\mapsto(\widetilde{M},\widetilde{g})$ is a conformal diffeomorphism between spaces of constant Riemannian curvature $C$ and $\widetilde{C}$. We can locally introduce co-ordinates $(x_1,\cdots,x_n)$ on $M$ and $(\widetilde{x}_1,\cdots,\widetilde{x}_n)$ on $\widetilde{M}$ which bring the metrics in so-called Riemannian form:
\begin{equation}
\label{eq:riem_form}
g = \frac{\dd x_1^{\,2}+\cdots+\dd x_n^{\,2}}{\left( 1 + \frac{C}{4}\sum(x_i)^2 \right)^2}
\qquad \textrm{and} \qquad
\widetilde{g} = \frac{\dd \widetilde{x}_1^{\,2}+\cdots+\dd \widetilde{x}_n^{\,2}}{\left( 1 + \frac{\widetilde{C}}{4}\sum(\widetilde{x}_i)^2 \right)^2}\,.
\end{equation}
It can easily be seen that the co-ordinate patches $x:M\rightarrow\mathbb{E}^n$
and $\widetilde{x}:\widetilde{M}\rightarrow\mathbb{E}^n$ are concircular, such that the diffeomorphism $\Psi$ is concircular if and only if its co-ordinate representation is a concircular diffeomorphism from $\mathbb{E}^n$ to $\mathbb{E}^n$. For $n > 2$ it easily follows from Liouville's theorem on conformal diffeomorphisms that this co-ordinate representation has to be a M\"obius transformation. For $n=2$ the same conclusion holds, but it is now a consequence of a result of Carath\'eodory \cite{caratheodory1937}. Of course, every diffeomorphism $\Psi$ the co-ordinate representation of which is a M\"obius transformation will be concircular.
\end{theorem_}

It can be remarked that the stereographic projection from part of the sphere $\mathrm{S}^n\subseteq\mathbb{E}^{n+1}$ to a hyperplane is such a co-ordinate system $x$ in which the spherical metric is represented in the form (\ref{eq:riem_form}).
A similar projection can be found for $\textrm{H}^n\subseteq\mathbb{E}^{n+1}_1$.

(Note that \cite{kuhnel1988}, cor. 10, for the case of spaces of constant sectional curvature and for $n=2$, is not correct: in this case there holds (ii) $\Rightarrow$ (i), but not conversely. The origin of the problem is a similar flaw in the sentence after \cite{kuhnel1988}, prop. 3.)


\section{Diffeomorphisms which Preserve Geodesic Hyperspheres.}
\label{ref:geodhyp}

We now ask whether Beltrami's theorem can likewise be adapted to the \textit{geodesic hyperspheres}, these being defined as the distance hyperspheres, \textit{i.e.}, the loci of all points which are separated a certain distance $r$ from a certain point $p$. This distance is to be measured by the Riemannian metric and this geodesic hypersphere will be denoted by $\mathscr{G}_p (r)$. For $r$ strictly positive and sufficiently small, this is a smooth hypersurface of the ambient Riemannian manifold.

We remark that several authors have already studied related questions, such as \cite{chenvanhecke,gray1973,grayvanhecke}: ``\textit{To which amount is the  Euclidean space characterised among all Riemannian  manifolds by the volume of its geodesic hyperspheres?}''

In relation with a possible adaption of Beltrami's theorem, we ask the similar, but simpler, question: ``\textit{To which amount is the Euclidean space characterised among all Riemannian  manifolds by its geodesic hyperspheres?}''

Before we can give an answer to this question in Theorem~\ref{thm:cosph} below, we will need two lemma's.

It should be mentioned that another approach towards this question has been suggested in \cite{blaschke_1954}, and a weaker version of Theorem~\ref{thm:cosph} can be found in \cite{dallavolta}.

\begin{lemma_}
\label{lem:cospher1}
Assume on a Riemannian manifold $M$ two metrics $g$ and $\widetilde{g}$ which share their geodesic hyperspheres are given. Consider, for a point $p\in M$, the one-parameter family of geodesic hyperspheres w.r.t.\ $(M,\widetilde{g})$ centred at $p$. 
Every such geodesic hypersphere $\widetilde{\mathscr{G}}_p(r)$ is also a geodesic hypersphere w.r.t.\ the metric $g$, although the centre $\beta(r)$ and the radius $\rho(r)$, as defined w.r.t. $(M,g)$, can be different. As such we have the relation $\widetilde{\mathscr{G}}_p(r)=\mathscr{G}_{\beta(r)}(\rho(r))$ for $r>0$ sufficiently small.

Now the extensions of the curve $\beta$ and the function $\rho$ which are obtained by setting $\beta(-r)=\beta(r)$ and $\rho(-r)=-\rho(r)$, are smooth in a neighbourhood of $0$.
\end{lemma_}

\begin{figure}
\begin{center}
\framebox{\includegraphics[width=.97\textwidth]{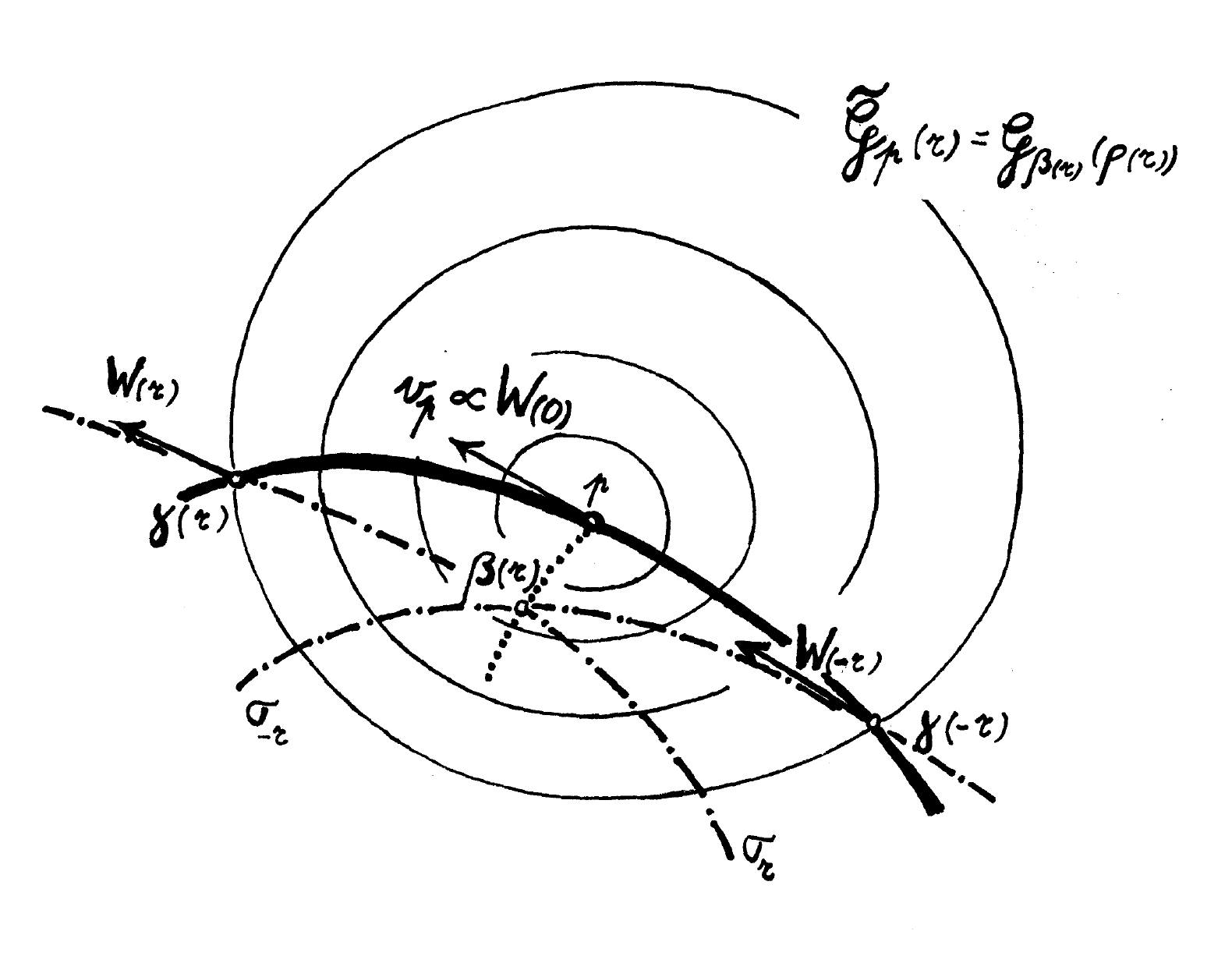}}
\caption{The geodesic hypersphere $\widetilde{\mathscr{G}}_p(r)$ of $(M,\widetilde{g})$ which is centred at $p$, is also a geodesic hypersphere of $(M,g)$. The centre $\beta(r)$ of this geodesic hypersphere, as determined w.r.t.\ $g$, lies on the intersection of the curves $\sigma_r$ and $\sigma_{-r}$.}
\label{fig:x1}
\end{center}
\end{figure}

\begin{figure}
\begin{center}
\framebox{\includegraphics[width=.97\textwidth]{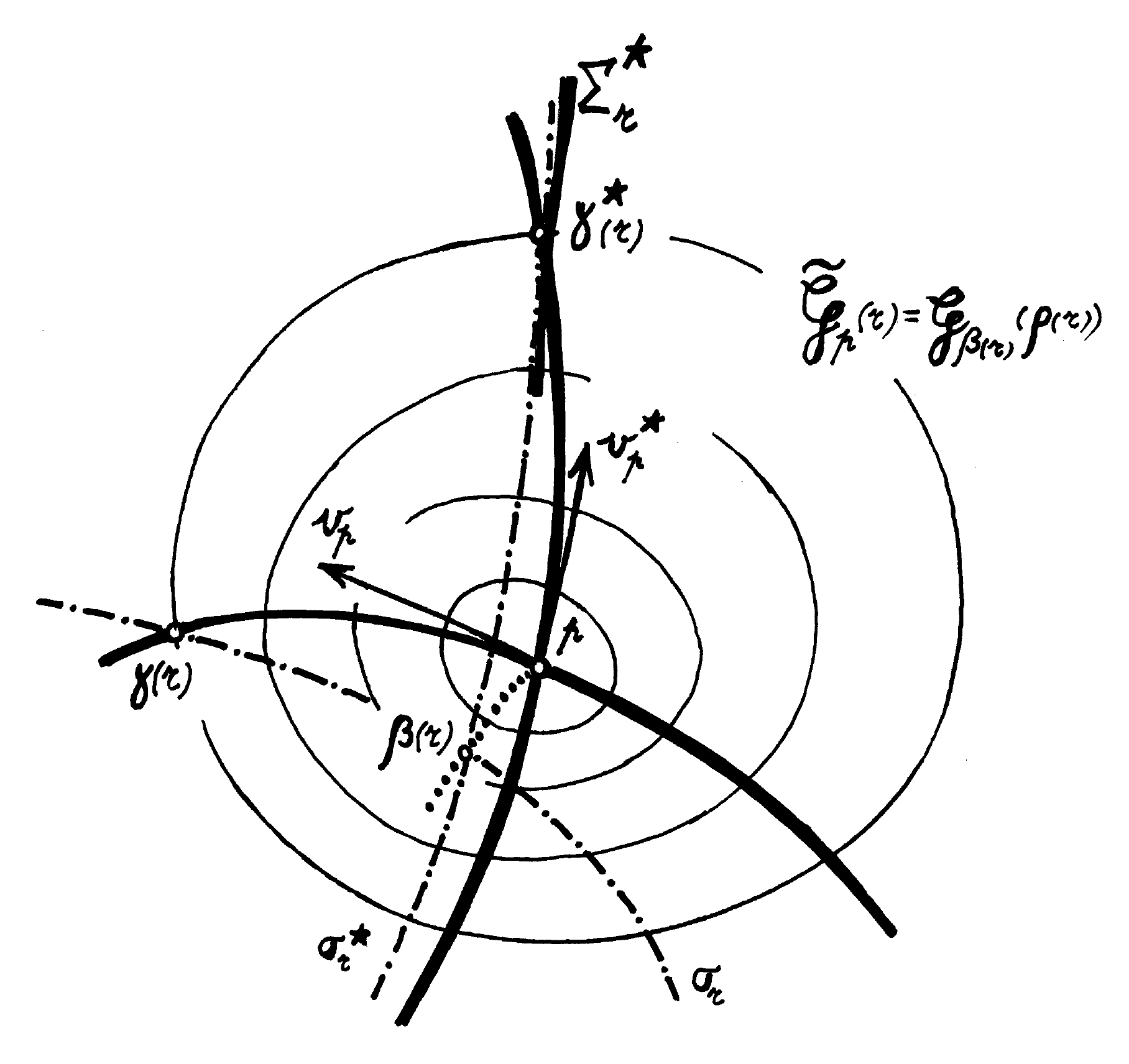}}
\caption{
The centre $\beta(r)$ of the geodesic hypersphere $\widetilde{\mathscr{G}}_p(r)$, viewed as a geodesic hypersphere of $(M,g)$, is uniquely determined as the point of intersection of the curve $\sigma_r$ and the hypersurface $\sigma_{r}^{\bigstar}$.}
\label{fig:x2}
\end{center}
\end{figure}

\begin{proof}
Let us first denote $\sigma$ for the unique operator satisfying $\widetilde{g} = \sigma \lrcorner\,g$. Choose an eigenvector $v_p\in\mathrm{T}_p M$ of $\sigma$ for which $\widetilde{g}(v_p,v_p)=1$, and let $\gamma$ be the $\widetilde{g}$--geodesic starting in $p$ with velocity $v_p$. Remark that the curve $\gamma$ and the hypersurface $\widetilde{\mathscr{G}}_p(r)$ cut each other orthogonally w.r.t.\ the metric $\widetilde{g}$ in the points $\gamma(\pm r)$. Consequently, if we define a vector field $W$ along $\gamma$ by
\[
W(r) = \frac{1}{\sqrt{\,\widetilde{g}(\gamma'(r),\sigma(\gamma'(r)))}\,}\,\sigma(\gamma'(r))
\]
then $W(r)$ stands orthogonally to this hypersurface
$\widetilde{\mathscr{G}}_p(r) = \mathscr{G}_{\beta(r)}(\rho(r))$ w.r.t.\ the metric $g$. Let us now denote, for every $r$ sufficiently small in absolute value, by $\sigma_r$ the $g$-geodesic with initial condition $\sigma_r'(0)=W(r)$. Remark that $\sigma_r(s)=\mathrm{exp}_{\gamma(r)}(s\,W(r))$ depends smoothly on $r$ and $s$ as a composition of smooth functions.

As is suggested in figure~\ref{fig:x1} on page~\pageref{fig:x1}, the point $\beta(r)$, for $r>0$, lies both on $\sigma_r$ and $\sigma_{-r}$.
To overcome the difficulty that $\beta(r)$ might even locally not be the unique point of intersection of these two curves, and to establish the smoothness of $\beta$, we will use a similar construction
which starts with a vector $v_p^{\bigstar}\in\mathrm{T}_p M$
for which $\widetilde{g}(v_p^{\bigstar},v_p^{\bigstar})=1$ and $\widetilde{g}(v_p^{\bigstar},v_p)=0$.
Define $\gamma^{\bigstar}$ as the $\widetilde{g}$--geodesic starting along the vector $v_p^{\bigstar}$.

Let $\Sigma^{\bigstar}_r$ stand for the hyperplane of $\mathrm{T}_{\gamma^{\bigstar}(r)}M$ which is the $g$--orthogonal complement of the vector which is obtained by $\widetilde{g}$--parallel transport of the vector $v_p\in\mathrm{T}_p M$ along $\gamma^{\bigstar}$. We define $\sigma_r^{\bigstar}$ as the hypersurface of $M$ which is spanned by all $g$--geodesics emanating from $\gamma^{\bigstar}(r)$ in a direction which is tangent to $\Sigma^{\bigstar}_r$. (See figure~\ref{fig:x2} on page~\pageref{fig:x2}.)

Now it will be seen that there exists a neighbourhood of the point $p$ such that for every $r$ sufficiently small in absolute value, $\beta(r)$ is the unique point of intersection of $\sigma_r$ and $\sigma_r^{\bigstar}$.

Choose a function $f_r:M\rightarrow \mathbb{R}$, smoothly depending on $r$, such that for every $r$ and every point $q$ in a neighbourhood of $p$ we have
\[
q \in \sigma_r^{\bigstar} \qquad \Leftrightarrow \qquad
f_r(q)=0\,,
\]
and such that $\dd f_r$ has no critical points in this neighbourhood.
The function $\rho$, as defined in the formulation of this lemma, obviously satisfies 
\begin{equation}
\label{eq:fsigmanul}
f_r(\sigma_r(-\rho(r)))=0
\end{equation}
for $r> 0$. Now define a function 
\[
F : \mathbb{R} \times \mathbb{R} \rightarrow \mathbb{R} : (x,y)\mapsto F(x,y)=f_x(\sigma_x(-y)).
\]
It should be remarked that $F_y(0,0)=\dd f_0(-\sigma'_0(0))=-\dd f_0(W(0))\propto -\dd f_0 (v_p)\neq 0$. Hence, according to the implicit function theorem, on some open interval around zero there uniquely exists a smooth, real-valued function $\rho$ for which
\[
F(x,\rho(x)) = 0 \qquad \textrm{and} \qquad\rho(0)=0.
\]
Because of (\ref{eq:fsigmanul}) this is an extension of the function $\rho$ which was originally defined. Then the curve $\beta$ can be smoothly extended left from zero by setting $\beta(r)=\sigma_r(-\rho(r))$.

It is clear from the construction that $\beta(-r)$ coincides with $\beta(r)$.
Since the $g$-geodesic distance between $\beta(\pm r)$ and $\gamma(\pm r)$ is $|\rho(\pm r)|$, it follows that $\rho(r)=-\rho(-r)$.
\end{proof}

The following lemma is related with results of M\"obius and  Liouville in the already interesting case when both $g$ and $\widetilde{g}$ are flat metrics (see, \textit{e.g.}, \cite{blair2000}, thm. 5.6; \cite{blaschke}, \S\,49).

\begin{lemma_}
\label{lem:cospher2}
Every cospherical diffeomorphism is conformal.
\end{lemma_}
\begin{proof}
We can assume that we are in the situation described in the previous lemma.

Now introduce an arbitrary co-ordinate system for which $p$ obtains co-ordinates $(0,\ldots,0)$. There holds
\begin{equation}
\label{eq:twoequations}
\left\{
\begin{array}{rcrcl}
\mathscr{G}_{\beta(r)}(\rho(r))&
\longleftrightarrow &
(\rho(r))^2 
&=&
\sum g_{ij\,(\beta(r))}\big(x_i- (\beta(r))_i \big)\,\big(x_j- (\beta(r))_j \big)
+\Order(|x|^3)
\,;\\
\widetilde{\mathscr{G}}_p(r) &
\longleftrightarrow &
r^2
&=& 
\sum \widetilde{g}_{ij\,(p)} x_i\, x_j
+\Order(|x|^3)
\,. \rule{0pt}{16pt}\\
\end{array}
\right.
\end{equation}
Define, for $y>0$ sufficiently small, $\tau(y)= (\rho(\sqrt{y}))^2$.
Because $\rho$ is an odd and smooth function, as was shown in the previous lemma~\ref{lem:cospher1}, this function $\tau$ is can be extended to a smooth function on an interval around $0$.
Since  the above equations (\ref{eq:twoequations}) are equivalent, we must have that for all $r$, for every point of $\mathscr{G}_{\beta(r)}(\rho(r))$, the following equation is satisfied:
\[
\tau\big( \sum \widetilde{g}_{ij\,(p)} x_i\, x_j
+\Order(|x|^3) \big)
\,=\,
\sum g_{ij\,(\beta(r))}\big(x_i- (\beta(r))_i \big)\,\big(x_j- (\beta(r))_j \big)
+\Order(|x|^3)
\,.
\]
Taking into account that $\beta'(0)=0$, this can be rewritten as
\begin{eqnarray*}
&&\sum \tau'(0)\,\widetilde{g}_{ij\,(p)} x_i\, x_j
\\
&&\quad=\quad 
\sum 
\Big(g_{ij\,(p)}+\frac{r^2}{2}\,(g_{ij}\circ\beta)''(0)  \Big)
\Big(x_i-\frac{r^2}{2}\,\beta_i''(0) \Big)\,
\Big(x_j-\frac{r^2}{2}\,\beta_j''(0) \Big)
+\Order(|x|^3)\,. 
\end{eqnarray*}
Interpreting $r:M\rightarrow \mathbb{R}^+$ as the geodesic distance to $p$ w.r.t.\ the metric $\widetilde{g}$, which has a co-ordinate expression already given in the second equation of (\ref{eq:twoequations}), the above equation should be identically satisfied.
From a comparison of the quadratic terms in $x$, we conclude that 
$\tau'(0)\,\widetilde{g}_{ij\,(p)}=g_{ij\,(p)}$.
\end{proof}

\begin{theorem_}[``\textit{Cospherical Beltrami theorem}'']
\label{thm:cosph}
If a cospherical diffeomorphism $\Psi$ between two Riemannian manifolds, one of which has constant Riemannian curvature, exists, then this diffeomorphism is conformal and both spaces have constant Riemannian curvature.

Moreover, for spaces of constant Riemannian curvature, the list of all cospherical diffeomorphisms coincides with the list of concircular diffeomorphisms, which has been given already.
\end{theorem_}

\begin{proof}
Let $g$ be a metric of constant Riemannian curvature $C$ on a manifold $M$ and $\widetilde{g}$ be another Riemannian metric on the same manifold $M$ such that the geodesic hyperspheres of both metrics coincide. A co-ordinate system on $M$ can be introduced for which $g$ takes the Riemannian form (\ref{eq:riem_form}). But then the co-ordinate system obviously realises a cospherical diffeomorphism between the space of constant Riemannian curvature and the parameter space $\mathbb{R}^n$ (with Euclidean metric).

For this reason we can simply assume that $g=\delta_{ij}$ is the Euclidean metric and, because of the previous lemma, that the second metric has the form $\widetilde{g}= \textrm{exp}(2\,\varphi)\,g$. We refer to \cite{kuhnel1988} for the standard formulae expressing the connection and the curvature of $\widetilde{g}$ in terms of $g$. We will use $\llangle\cdot,\cdot\rrangle$ as alternative notation for this Euclidean metric $g$ and denote $\|\cdot\|$ for the corresponding norm. Differential invariants w.r.t. the metric $\widetilde{g}$ will always bear a tilde in their notation.

We will keep the notation of the previous proofs where now $M$ becomes $\mathbb{R}^n$.

It can be seen that the curve $\gamma:\mathbb{R}\rightarrow \mathbb{R}^n  : t\mapsto \gamma(t)$ is a geodesic w.r.t.\ $\widetilde{g}$ if and only if
\[
\gamma '' = \llangle \gamma' , \gamma' \rrangle \textrm{grad} \varphi - 2\,\llangle  \textrm{grad}\varphi ,\gamma'\rrangle \gamma'\,.
\]
(Here $\gamma''$ is just the second derivative of $\gamma$ as a curve in $\mathbb{R}^n$\!.)
If $v\in\mathrm{T}_p M$ with $\| v \| =1$, then the corresponding $\widetilde{g}$--geodesic satisfies consequently
\begin{eqnarray*}
\gamma_v(t) &=& p+ t\,v + \frac{t^2}{2} \left\lgroup \textrm{grad}\varphi-2\llangle \textrm{grad}\varphi,v\rrangle v \right\rgroup \\
&&+\frac{t^3}{6}\left\lgroup
\textrm{Hs}_{\varphi}(v) - 2\,\textrm{Hess}_{\varphi}(v,v)\,v - 2\, \|\textrm{grad}\varphi\|^2\, v\right.\\
&&\qquad\qquad\qquad\qquad\qquad\left. + 8\, \llangle \textrm{grad}\varphi ,v\rrangle^2 v - 4\, \llangle \textrm{grad}\varphi,v \rrangle \textrm{grad}\varphi
\right\rgroup
+
\Order(t^3)\,.
\end{eqnarray*}
Here it has been suppressed in the notation that $\textrm{grad}\varphi$, $\textrm{Hs}_{\varphi}$, $\cdots$ have to be evaluated at the point $p$.
Furthermore we remark that $ \gamma_v(t)= \widetilde{\mathrm{exp}}_p(t\,v)$. If, for some fixed $t$, the vector $v$ runs through all $v\in\mathrm{T}_p M$ with $\| v \| =1$, \textit{i.e.}, with $\sqrt{\widetilde{g}(v,v)}=\mathrm{exp}(\varphi(p))$, the point $\gamma_v(t)$ will describe the geodesic hypersphere 
\[
\widetilde{\mathscr{G}}_p (t\, \mathrm{exp}(\varphi(p)))= \mathscr{G}_{\beta(t\, \mathrm{exp}(\varphi(p)))}(\rho(t\, \mathrm{exp}(\varphi(p))))\,,
\]
\textit{i.e.}, a Euclidean hypersphere of a certain radius and center
\[
{\beta(t\, \mathrm{exp}(\varphi(p)))} = p+ \frac{t^2}{2} b_2 + \frac{t^4}{24} b_4 + \cdots\,.
\]
It has already been shown that no odd terms occur in the right-hand side of the above equation.

This means that, for every $t$, the number
\[
\llangle \left(\gamma_v(t) -p- \frac{t^2}{2}b_2-\frac{t^4}{24}b_4 -\cdots\right)
\, , \,
\left(\gamma_v(t) -p- \frac{t^2}{2}b_2-\frac{t^4}{24}b_4 - \cdots \right)\rrangle
\]
is independent of the choice of $v\in \mathrm{T}_p M$ with $\|v\|=1$. Hence if we develop the above expression in powers of $t$, all coefficients should be independent of the choice of  $v\in \mathrm{T}_p M$ with $\|v\|=1$.
This condition is automatically satisfied for the coefficient of $t^2$.
For the coefficient of $t^3$ this condition means that $b_2=-(\mathrm{grad}\varphi)_{(p)}$.
The coefficient of $t^4$ is given by 
\[
\frac{1}{3}\left\lgroup-\mathrm{Hess}_{\varphi}(v,v)+\,\|\mathrm{grad}\varphi\|^2 
+ \llangle \mathrm{grad}\varphi,v\rrangle^2\right\rgroup.
\]
The fact that this coefficient is independent of $v\in \mathrm{T}_p M$ (subject to 
$\|v\|=1$), and this no matter how the point $p$ has been chosen, means that 
\begin{equation}
\label{eq:phi}
\mathrm{Hess}_{\varphi} = \mu\, g + \dd \varphi \otimes \dd\varphi
\end{equation}
for some function $\mu$. If the function $h$ is introduced by the relation $h^{-2}=\mathrm{exp}(2\varphi)$, the above condition (\ref{eq:phi}) precisely means that $\mathrm{Hess}_{h} = \frac{\Delta h}{n} \,g$. A comparison with \cite{kuhnel1988}, prop. 8, shows that the diffeomorphism is also concircular,
which completes the proof.
\end{proof}


\section{Diffeomorphisms which Preserve Minimal Surfaces.}

Before we can determine all cominimal diffeomorphisms between three-dimensional Riemannian manifolds, we need two lemma's.

\begin{lemma_}
\label{lem:comin1}
Let $M$ be a surface of a three-dimensional Riemannian manifold $\overline{M}$, and assume that the locus of points where the mean curvature of $M\subseteq\overline{M}$ vanishes is a smooth curve $\gamma$.

The image of this curve under a cominimal diffeomorphism from $\overline{M}$ to another three-dimensional Riemannian manifold is the locus of points where the mean curvature of the image surface vanishes.
\end{lemma_}
\begin{proof}
Consider the first-order surface band (``Fl\"achenstreifen'') $\mathfrak{B}$ determined by the curve $\gamma$ and the tangent spaces of $M$ along this curve. This band $\mathfrak{B}$ provides Cauchy initial data for the minimal surface equation, and, due to the Cauchy--Kowalewski theorem, locally determines a unique minimal surface $M_{\ast}$ in $\overline{M}$ (see \cite{leichtweiss1956}).

The fact that the mean curvature of $M$ 
vanishes along the curve $\gamma$
precisely means that $M$ and $M_{\ast}$ have higher-order contact along the curve $\gamma$.
 
This notion of higher-order contact is, of course, preserved under the cominimal diffeomorphism $\Psi$. Since $\Psi(M_{\ast})$ is a minimal surface in the second Riemannian manifold which has higher-order contact with $\Psi(M)$ along $\Psi(\gamma)$, the result follows.
\end{proof}

\begin{lemma_}
\label{lem:comin2}
Let $(\overline{M},\overline{g})$ be a three-dimensional Riemannian manifold. Assume for a point $p\in\overline{M}$, two orthonormal vectors in $\mathrm{T}_p \overline{M}$ are given. For every real number $\ell$, there exists a surface $M\subseteq\overline{M}$ which passes through $p$, which has the two given vectors as principal directions at $p$, and $\pm\ell$ as principal curvatures.
\end{lemma_}
\begin{proof}
It is clear that a surface $M_{p}$ in the Euclidean space 
$(\mathrm{T}_p \overline{M},\overline{g}_{p})$ exists which satisfies all requirements.
Now define $M$ as the image of $M_{p}$ under the Riemann exponential diffeomorphism of $(\overline{M},\overline{g})$. 
\end{proof}

\begin{theorem_}[``\textit{Cominimal Beltrami theorem}'']
\label{thm:comin}
There do not exist cominimal diffeomorphisms between three-dimensional Riemannian manifolds, except for homotheties.
\end{theorem_}

The above theorem can also be derived from Thm.\ II of \cite{lapaz} in case one of the spaces is Euclidean.

\begin{proof}
Let us start by explaining the notation. We assume that $\overline{g}$ and $\widetilde{\overline{g}}$ are two metrics on a three-dimensional manifold, which will be denoted by $\overline{M}$ or $\widetilde{\overline{M}}$ depending on which metric is involved. 

A surface in this space will usually be denoted by $M$ when it is regarded as submanifold of $(\overline{M},\overline{g})$, and as $\widetilde{M}$ when regarded as submanifold of $(\widetilde{\overline{M}},\widetilde{\overline{g}})$. We generally follow the notational rule that objects without resp. with a $\overline{\textrm{bar}}$ refer to the induced resp.\ the ambient geometry, and without resp. with a $\widetilde{\textrm{tilde}}$ to the metric $\overline{g}$ resp.\ $\widetilde{\overline{g}}$. 

For instance, for vector fields $V,W\in\mathfrak{X}(M)$ we will write the Gauss equation for the surface $M\subseteq\overline{M}$ as
\[
\overline{\nabla}_V W = \nabla_V W + \II(V,W)\,N\,;
\]
if we study the surface $\widetilde{M}\subseteq\widetilde{\overline{M}}$ the corresponding equation is
\[
\widetilde{\overline{\nabla}}_V W = \widetilde{\nabla}_V W + \widetilde{\II}(V,W)\,\widetilde{N}\,.
\]
We will decompose the unit normal vector field 
$\widetilde{N}$ (w.r.t.\ $\widetilde{\overline{g}}$\,) along the surface 
$\widetilde{M}$ in $\widetilde{\overline{M}}$ in a part which is proportional to the normal $N$ (w.r.t.\ $\overline{g}$) and a part which is tangent to the surface:
\[
\widetilde{N} = \varphi \,N + \widetilde{N}^{\tangentpart}\,.
\]
We also introduce the difference tensor $\overline{X}$ on $\overline{M}$  between the two connections by
\[
\overline{X}(V,W) = \widetilde{\overline{\nabla}}_V W -\overline{\nabla}_V W \,,
\]
for every $V,W\in\mathfrak{X}(\overline{M})$. With this notation, there holds
\begin{eqnarray*}
\widetilde{A}(V) 
&=& -\widetilde{\overline{\nabla}}_V \widetilde{N}
= -\overline{\nabla}_V \widetilde{N} - \overline{X}(V,\widetilde{N})\\
&=& -V[\varphi] N + \varphi\,A(V) - \nabla_V(\widetilde{N}^{\tangentpart})
-\II(V,\widetilde{N}^{\tangentpart})\,N
-\overline{X}(V,\widetilde{N})
\end{eqnarray*}
for every $V\in\mathfrak{X}(M)$.
If we consider the tangent part of both sides of this equation,  we obtain the following relation between the shape operators with respect to the different ambient metrics:
\[
\widetilde{A}(V) - \varphi\, A(V) = - \nabla_V (\widetilde{N}^{\tangentpart}) - \left(\overline{X}(V,\widetilde{N})\right)^{\tangentpart}\,,
\]
in which $v^{\tangentpart}$ stands for the tangent part of a vector $v$ w.r.t. the metric $\overline{g}$. By taking the trace there results
\begin{equation}
\label{eq:HHtilde}
2\,\widetilde{H} - 2\,\varphi\, H =
-\,\mathrm{div}(\widetilde{N}^{\tangentpart}) 
- \mathrm{trace}\left\{V \mapsto \left(\overline{X}(V,\widetilde{N})\right)^{\tangentpart} \right\}\,.
\end{equation}
Here the divergence and the trace are taken on $(M,g)$.
All of the above equations are valid for every surface in every three-dimensional Riemannian manifold endowed with two metrics $\overline{g}$ and $\widetilde{\overline{g}}$.
In this case there uniquely exists an operator $\sigma:\mathfrak{X}(\overline{M})\rightarrow\mathfrak{X}(\overline{M})$, symmetric w.r.t.\ $\overline{g}$, such that 
$\widetilde{\overline{g}} = \sigma \lrcorner\,\overline{g}$.
Let $s_1$, $s_2$ and $s_3$ be the eigenvectors of $\sigma$ (which have been normalised w.r.t.\ $\overline{g}$, and hence consitute a $\overline{g}$-orthonormal frame field on $\overline{M}$), and $\sigma_1$, $\sigma_2$ and $\sigma_3$
the corresponding eigenvalues. The latter are strictly positive since both metrics are positive-definite.

For a surface $M\subseteq\overline{M}$, we will denote $e_1$, $e_2$ for the principal directions and $\lambda_1$, $\lambda_2$ for the principal curvatures (w.r.t. $(\overline{M},\overline{g})$ as ambient space).

Now assume that these metrics share their minimal surfaces. 
Choose an arbitrary point $p\in \overline{M}$. It has been shown in lemma~\ref{lem:comin2} that a surface $M\subseteq\overline{M}$ exists which passes through $p$ with prescribed tangent plane $\mathrm{span}\{(s_1)_{(p)},(s_2)_{(p)}\}$, prescribed principal directions $(e_1)_{(p)}=(s_1)_{(p)}$ and $(e_2)_{(p)}=(s_2)_{(p)}$, and prescibed mean curvature $H_{(p)}=0$, and that the number 
$(\lambda_1)_{(p)}$ can additionally be chosen freely.

The unit normal $N$ of such a surface $M\subseteq\overline{M}$ will be chosen in such a way that $N_{(p)}=(s_3)_{(p)}$, and the unit normal $\widetilde{N}$ can be chosen as
\[
\widetilde{N} = \frac{1}{\sqrt{\overline{g}(N,\sigma^{-1}(N)) }} \,
\sigma^{-1}(N) \,.
\]
For such a surface, the following equations hold true (in which all sums run over $i=1,2$):
\begin{eqnarray*}
&\!\!\mathrm{div}&\!\!(\widetilde{N}^{\tangentpart})
= \sum\, g\left( \nabla_{e_i} \left(
\frac{1}{\sqrt{\overline{g}(N,\sigma^{-1}(N)) }} 
\left(\sigma^{-1}(N) - \overline{g}(\sigma^{-1}(N),N)N \right)
\right) , e_i \right)\\
&=& 
\sum\, \Bigg\{ 
e_i\left[\frac{1}{\sqrt{\overline{g}(N,\sigma^{-1}(N)) }}\right]
\overline{g}(\sigma^{-1}(N),e_i)
\\
&&
\quad\qquad
+\frac{1}{\sqrt{\overline{g}(N,\sigma^{-1}(N)) }}
\,\overline{g}\left(\overline{\nabla}_{e_i}\left(\sigma^{-1}(N)\right) ,e_i\right)\\
&&
\quad\qquad-\sqrt{\overline{g}(N,\sigma^{-1}(N)) } 
\,\overline{g}\left(\overline{\nabla}_{e_i} N ,e_i\right)
\Bigg\}\\
&=&  \left(\big(\sigma^{-1}(N)\big)^{\tangentpart} \right)\!\!\left[\frac{1}{\sqrt{\overline{g}(N,\sigma^{-1}(N)) }}\right]\\
&&
+\sum\, \frac{1}{\sqrt{\overline{g}(N,\sigma^{-1}(N)) }} 
\Bigg\{
-\overline{g}\left(
\sigma^{-1}\left((\overline{\nabla}_{e_i}\sigma)(\sigma^{-1}(N)) \rule{0pt}{12pt} \right) 
, e_i
\right)
+\overline{g}\left(\sigma^{-1}(\overline{\nabla}_{e_i}N),e_i \right)    \Bigg\}\\
&&
-\sum\,\sqrt{\overline{g}(N,\sigma^{-1}(N)) } 
\,\overline{g}\left(\overline{\nabla}_{e_i} N ,e_i\right)
\,.\rule{0pt}{18pt}
\end{eqnarray*}
If we evaluate both sides of this equation at the point $p$, where the relations 
$(e_i)_{(p)}=(s_i)_{(p)}$ (for $i=1,2$) and $N_{(p)}=(s_3)_{(p)}$ are satisfied, we obtain
\begin{eqnarray}
\label{eq:divNrhs}
\left(\mathrm{div}(\widetilde{N}^{\tangentpart})\right)_{\!(p)}
&=&
\left\lgroup\rule{0pt}{18pt}\right. -\sum\, 
\sqrt{\sigma_3}\,
\Bigg\{
\overline{g}\left(
\sigma^{-1}\left((\overline{\nabla}_{s_i}\sigma)(\sigma^{-1}(s_3)) \rule{0pt}{12pt} \right) 
, s_i
\right)\\
\nonumber
&&\qquad\qquad\qquad
+\overline{g}\left(\sigma^{-1}(A(e_i)),e_i \right)    \Bigg\}
+ \sum\,\frac{1}{\sqrt{\sigma_3}}\, \overline{g}(A(e_i),e_i) \left.\rule{0pt}{18pt}\right\rgroup_{\!\!(p)}\\
\nonumber
&=&
\left\lgroup\rule{0pt}{18pt}\right. 
-\sum\, 
\frac{1}{\sigma_i\sqrt{\sigma_3}}\,
\overline{g}\left(\left(\overline{\nabla}_{s_i} \sigma\right)(s_3),s_i\rule{0pt}{12pt}  \right)
+\sum\,\sqrt{\sigma_3} \left(\frac{-1}{\sigma_i}+\frac{1}{\sigma_3}\right)\lambda_i\left.\rule{0pt}{18pt}\right\rgroup_{\!\!(p)}\,.
\end{eqnarray}
The only restrictions which we have imposed on the surface with vanishing mean curvature in $p$ was its tangent plane at $p$ as well as its principal directions $(e_1)_{(p)}$ and $(e_2)_{(p)}$. Consider a family of  surfaces which satisfy all these conditions but for which $(\lambda_1)_{(p)}$ takes different values. Due to lemma~\ref{lem:comin1}, the left-hand side of (\ref{eq:HHtilde}), when evaluated at $p$, vanishes for any of these surfaces. Moreover, the second term in the right-hand side of this equation, when evaluated at $p$,  will not depend on the actually chosen surface because it is algebraic in $\widetilde{N}$. Consequently, $\left(\mathrm{div}(\widetilde{N}^{\tangentpart})\right)_{(p)}$ will not depend on the chosen surface either. If we look at equation
(\ref{eq:divNrhs}), we see that the only term which can depend on the surface is the last one.

Therefore we must have that
\[
\left\lgroup\rule{0pt}{18pt}\right. 
\sum\,\sqrt{\sigma_3} \left(\frac{-1}{\sigma_i}+\frac{1}{\sigma_3}\right)\lambda_i\left.\rule{0pt}{18pt}\right\rgroup_{\!\!(p)}
=
\left\lgroup\rule{0pt}{18pt}\right. 
\lambda_1\,\sqrt{\sigma_3} \left(\frac{1}{\sigma_2}-\frac{1}{\sigma_1}\right)
\left.\rule{0pt}{18pt}\right\rgroup_{\!\!(p)}
\]
is independent of $(\lambda_1)_{(p)}$. This occurs only if $(\sigma_1)_{(p)}=(\sigma_2)_{(p)}$.

Hence, $\sigma$ is a multiple of the identity, \textit{i.e.}, the two metrics $\overline{g}$ and $\widetilde{\overline{g}}$ are conformal. This implies $\widetilde{N}^{\tangentpart}=0$,
and from eq. (\ref{eq:HHtilde}) can now be seen that, for all $W\in\mathfrak{X}(\overline{M})$,
\[
0 = \mathrm{trace}\left\{V \mapsto \left(\overline{X}(V,W)\right)^{\tangentpart} \right\}
\]
holds. Taking the expression for the behaviour of the Levi-Civita connection under a conformal transformation into account, we find that the conformal factor has to be constant.
\end{proof}


\section{The infinitesimal Beltrami Theorem.}

\subsection{Definition of an Infinitesimally Cogeodesical Deformation.}
\label{sec:def_infproj}
Consider a one-parameter family of Riemannian metrics on an $n$-dimensional manifold $M$:
\[
g^{(t)}=g+t\,\delta g+\Order(t^2) \qquad\qquad \textrm{(\,$t\in\left]-\varepsilon,\varepsilon\right[$\,)}\,.
\]
Since we can construct the usual Riemannian invariants with respect to any of the metrics $g^{(t)}$, we will distinguish them by adding an index $(t)$ in the notation. However for $t=0$ this index will usually be omitted.

\begin{figure}
\begin{center}
\framebox{\includegraphics[width=.97\textwidth]{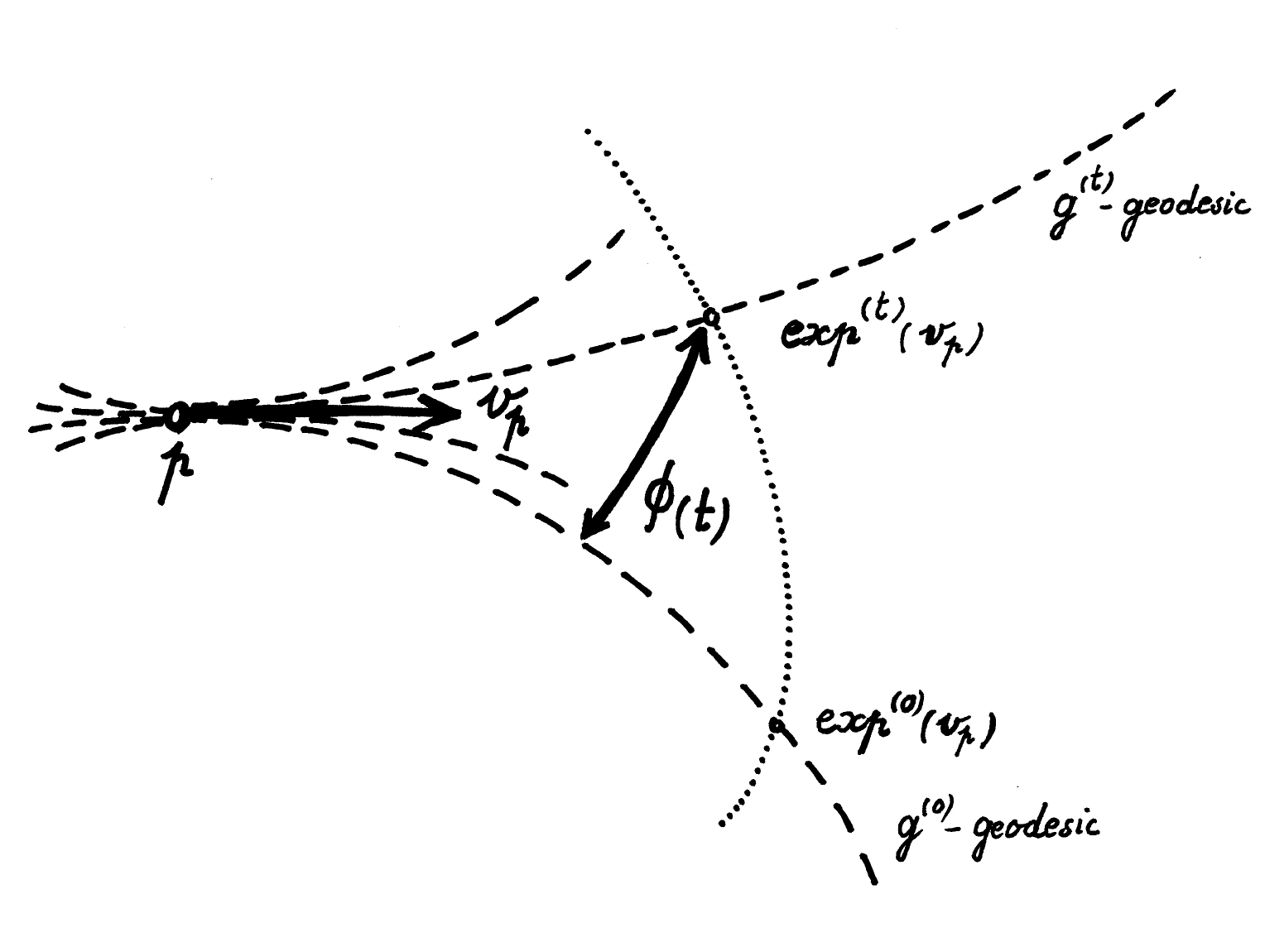}}
\caption{The distance $\phi(t)$ between the point $\exp^{(t)}(v_p)$ and the $g$-geodesic which starts in $p$ with velocity $v_p$.}
\label{fig:geod}
\end{center}
\end{figure}

\begin{definition_}
\label{def:inf_proj_def}
In the above situation, choose $p\in M$ and $v_p\in\mathrm{T}_p M$ arbitrarily and let $\phi(t)$ stand for the distance (measured with respect to the initial metric $g$) between the point $\exp^{(t)}(v_p)$ and the $g$-geodesic which starts in $p$ with velocity $v_p$ (see figure \ref{fig:geod} on p.\ \pageref{fig:geod}). If for all $p\in M$ and all $v_p\in\mathrm{T}_p M$ this distance satisfies
\[
\phi(t)=0\qquad \textrm{resp.} \qquad
\phi(t)=0+\Order(t^2)\,,
\]
then the family $\{\,g^{(t)}\,\}$ is called a \textit{cogeodesical deformation} resp.\ an \textit{infinitesimally cogeodesical deformation} of $(M,g)$.
\end{definition_}

Thus a cogeodesical deformation is a one-parameter family of geodesically equivalent metrics and an infinitesimally cogeodesical deformation is a one-parameter family of metrics which are, in a geometric sence, ``up to first-order approximation'' geodesically equivalent.

\subsection{Characterisation of Infinitesimally Cogeodesical Deformations by means of the Variation Tensor of the Connection.}
\label{sec:lem_infproj}

The \textit{variation tensor of the connection} is denoted by $X$ and defined by $X(V,W)=\deltatnought \nabla_V^{(t)} W$ for all $V, W \in \mathfrak{X}(M)$.
Furthermore, let $\sigma$ stand for the operator for which $\delta g = \sigma \lrcorner\,g$ holds true.

\begin{lemma_}
\label{lem:beltrami}
A family of metrics $\left\{g^{(t)}\right\}$ of an $n$-dimensional Riemannian manifold $(M,g)$ (with $g^{(0)}=g$) is an infinitesimally cogeodesical deformation if and only if the following equation is satisfied for all $V,W\in\mathfrak{X}(M)$:
\begin{equation}
\label{eq:X}
X(V,W) = \frac{1}{2\,(n+1)} V[\,\mathrm{trace}\, \sigma\,]\, W
+\frac{1}{2\,(n+1)} W[\,\mathrm{trace}\,\sigma\,]\, V\,.
\end{equation}
\end{lemma_}

\begin{proof}
Choose a point $p\in M$ and a vector $v_p\in \mathrm{T}_p M$ arbitrary. We consider the one-parameter family of curves $\{\,\gamma_t\,\}$ as on the figure,
\[
\gamma_t(s)=\exp^{(t)}(s\,v_p)\,.
\]
This family of curves can be regarded as a deformation of the curve $\gamma:=\gamma_0$, and the corresponding deformation vector field will be denoted by $V$:
\[
V_{(s)} = \deltatnought \gamma_t(s)\,.
\]
This is a $\mathrm{T}M$--valued vector field along $\gamma$ which is completely determined by the requirement
\begin{equation}
\label{eq:nablatgamma}
0= \nabla^{(t)}_{\gamma_{t}'(s)}\gamma_{t}'\,,
\end{equation}
and the initial conditions $\gamma_{t}(0)=p$ and $\gamma_{t}'(0)=v_p$.
Let us decompose this vector field as
\begin{equation}
\label{eq:97_2}
V = \varphi\,\gamma' + V^{\perp}
\end{equation}
where $g(\gamma',V^{\perp})=0$ and $\varphi$ is a function of the parameter $s$ of the curve $\gamma$.
It is straightforward to see that $\left\{g^{(t)}\right\}$ is an infinitesimally cogeodesical deformation if and only if $V^{\perp}=0$ for all choices of $p$ and $v_p$.

Consider now the following two vector fields $V$ and $W$ along the two-parameter mapping $\mathbf{x}:\mathbb{R}^2\rightarrow M : (t,s)\mapsto \gamma_t(s)$:
\[
V = \dd\mathbf{x} (\partial_t) \qquad \mathrm{and} \qquad W = \dd\mathbf{x}(\partial_s)\,.
\]
This extends a previous definition of $V$. Because these two vector fields originate from a two-parameter mapping, the Lie bracket $[V,W]$ vanishes.
Furthermore, we claim that 
\begin{equation}
\label{eq:claim}
\nabla_V \left(\nabla_W W\right) = -X(W,W) 
\end{equation}
holds along $\gamma$.
To see this, choose $s$ fixed, and define
\[
\eta:\mathbb{R}\rightarrow M : t \mapsto \gamma_t(s) = \exp^{(t)}(s\,v_p)\,. 
\]
For instance, for $s=1$, this becomes the dotted line 
($\cdot\!\cdot\!\cdot\!\cdot\!\cdot\!\cdot\!\cdot$) in the figure. We can rewrite equation (\ref{eq:nablatgamma}) as
\[
0= \nabla_{W_{\eta(t)}}^{(t)} W \,.
\]
Let now, for a vector $v\in\mathrm{T}_{\eta(t)}M$, the notation $\left.\left\{v\right\}\right\downarrow_{\eta(0)}$ stand for the vector of $\mathrm{T}_{\eta(0)}M$ which is obtained from $v$ by parallel transport of $v$ w.r.t. the connection $\nabla$ along the curve $\eta$ back to the initial point $\eta(0)$. Then there also holds
\[
0= \left.\left\{\nabla_{W_{\eta(t)}}^{(t)} W \right\}\right\downarrow_{\eta(0)}
\]
and consequently we have
\begin{eqnarray*}
0
&=& 
\deltatnought
\left.\left\{\nabla_{W_{\eta(t)}}^{(t)} W \right\}\right\downarrow_{\eta(0)}\\
&=&
\deltatnought
\left.\left\{\nabla_{W_{\eta(0)}}^{(t)} W \right\}\right\downarrow_{\eta(0)}
+
\deltatnought
\left.\left\{\nabla_{W_{\eta(t)}}^{(0)} W \right\}\right\downarrow_{\eta(0)}\\
&=&
X(W_{(\gamma(s))},W_{(\gamma(s))}) + 
\nabla_{V_{(\gamma(s))}} \left(\nabla_W W\right) \,,
\end{eqnarray*}
which proves our claim (\ref{eq:claim}).
This implies the following equations (valid along $\gamma$), where the last step comes from the decomposition (\ref{eq:97_2}), and where a prime denotes covariant differentiation along the curve $\gamma$ with respect to the connection $\nabla$:
\begin{eqnarray*}
-X(\gamma',\gamma') 
&=& \nabla_V (\nabla_W W) = \nabla_W (\nabla_V W) - \mathrm{R}(V,W)W\\
&=& \nabla_W (\nabla_W V) - \mathrm{R}(V,W)W\\
&=& \nabla_{\gamma'} (\nabla_{\gamma'} V) - \mathrm{R}(V,\gamma')\gamma'\\
&=& \varphi''\, {\gamma'} +2\,\varphi'\, \gamma'' + \varphi\, \gamma''' + (V^{\perp})'' + \mathrm{R}(\gamma', V^{\perp}) \gamma'\,.
\end{eqnarray*}
The second and the third term vanish because $\gamma$ is a geodesic. We obtain the following equation, which will be needed later:
\begin{equation}
\label{eq:XTT}
X(\gamma',\gamma') = 
\underbrace{\raisebox{-8pt}{\rule{0pt}{8pt}}
-\phi''\,\gamma'}_{\textrm{tangent to $\gamma\phantom{'}$}}
\underbrace{\raisebox{-8pt}{\rule{0pt}{8pt}}
- \left( (V^{\perp})'' +  \mathrm{R}(\gamma', V^{\perp}) \gamma'\right)}_{\textrm{normal to $\gamma'$}}
\end{equation}

Using Koszul's formula, the following relation between the tensors $X$ and $\sigma$ can be found:
\begin{equation}
\label{eq:Xsigma}
2\,g(X(V,W),Z) = g((\nabla_V \sigma)W,Z)+g((\nabla_W \sigma)V,Z)-g((\nabla_Z \sigma)V,W)
\,.
\end{equation}
As an immediate consequence of this formula we obtain
\begin{equation}
\label{eq:tracedeltag}
\textrm{trace}X_V := 
\textrm{trace}\{\,W \mapsto X(V,W)\,\} = \frac{1}{2} V[\,\textrm{trace}\,\sigma\,]\,.
\end{equation}

We are now in a position to complete the proof of the lemma. Assume first that 
$\left\{g^{(t)}\right\}$ is an infinitesimally cogeodesical deformation, such that $V^{\perp}=0$. It follows from equation (\ref{eq:XTT}) that $X(\gamma',\gamma')\propto \gamma'$. A polarisation identity for the tensor $X$ then implies that there exist functions $\alpha$ and $\beta$ for which the equation
\[
X(V,W) = \alpha(V,W)\,W + \beta(V,W)\,V
\]
is satisfied. The $\mathbb{R}$-linearity of $X$ immediately implies that $\alpha(V,W)$ is independent of $W$ and is linear in $V$, and similar for $\beta$. Taking the symmetry of $X$ into account, we find that
\[
X(V,W) = \nu(V)\,W + \nu(W)\,V
\]
for a certain one-form $\nu$. A combination of this equation with (\ref{eq:tracedeltag}) then easily results in the final formula (\ref{eq:X}) for the tensor $X$.

Assume now conversely that this formula (\ref{eq:X}) is satisfied. Choose a vector $v_p\in \mathrm{T}M$ and let $\gamma$ be as above the geodesic starting at $p$ in direction $v_p$. It follows from (\ref{eq:X}) that $X(\gamma',\gamma')\propto \gamma'$, and hence equation (\ref{eq:XTT}) gives that $V^{\perp}$ is a Jacobi vector field along $\gamma$. 
From $V_{(0)}=\deltatnought \gamma_t(0)=\deltatnought p=0$ obviously 
follows $(V^{\perp})_{(0)}=0$. 
Taking into account that $\gamma$ is a geodesic, we have 
\[
\varphi'(0)\, \gamma'(0)+ (V^{\perp})'_{(0)}=V'_{(0)} = \left(\rule{0pt}{13pt}\nabla_{\partial_t}\left(\dd\mathbf{x}(\partial_s) \right)\right)_{(s,t)=(0,0)}= 0\,,
\]
because $\dd\mathbf{x}(\left.\partial_s\right\vert_{s=0})=\gamma_t'(0)=v_p$ independently of $t$. Consequently there holds $(V^{\perp})'_{(0)}=0$. It now follows that $V^{\perp}=0$. Since this holds for every $v_p$, $\left\{g^{(t)}\right\}$ is an infinitesimally cogeodesical deformation.
\end{proof}

\subsection{The Infinitesimal Beltrami Theorem for Dimension $n=2$.}
\label{sec:two_dim}

\begin{theorem_}[``\textit{Infinitesimal Beltrami theorem}'']
\label{thm:beltrami}
\label{thm:inf_beltrami}
The variation of the Gaussian curvature of a two-dimensional space of constant Gaussian curvature under an infinitesimally cogeodesical deformation is constant.
\end{theorem_}

\begin{proof}
\noindent\textbf{(i.)}---Assume first that the curvature of the initial metric $g$ on the two-dimensional manifold $M$ vanishes, such that we can find a co-ordinate system on $M$ for which
\[
g^{(t)} =
\left\lgroup
\begin{array}{cc}
E^{(t)} & F^{(t)}\\
F^{(t)} & G^{(t)}
\end{array}
\right\rgroup
\quad
\textrm{where}
\quad
g^{(0)} =
\left\lgroup
\begin{array}{cc}
1 & 0\\
0 & 1
\end{array}
\right\rgroup\,.
\]
Let us also write $\Gamma_{ij}^{k\,(t)}$ for the Christoffel symbols w.r.t. the metric $g^{(t)}$. For the variation tensor of the connection we have the co-ordinate expression $X_{ij}^k = \deltatnought \Gamma_{ij}^{k\,(t)}$. It immediately follows from formula (\ref{eq:X}) that
\begin{equation}
\label{eq:XX}
X_{11}^{1} = 2\, X_{12}^{2}\,; \qquad
X_{11}^{2} = 0\,; \qquad
X_{22}^{1} = 0\,; \qquad
X_{22}^{2} = 2\, X_{12}^{1}\,. 
\end{equation}
Now for every value of $t$, the Gaussian curvature $K^{(t)}$ of the two-dimensional Riemannian manifold $(M,g^{(t)})$ can be computed from the Christoffel symbols with help of the well-known equations
\begin{equation}
\label{eq:Kgij}
\!\!\!\!\!\!\!\!\left\{
\begin{array}{rcl}
K^{(t)}\,g_{11}^{(t)} &=&
\partial_2 \Gamma_{11}^{2\,(t)}
- \partial_1 \Gamma_{12}^{2\,(t)}
-\Gamma_{11}^{2\,(t)}\,\Gamma_{12}^{1\,(t)}
+\Gamma_{11}^{1\,(t)}\,\Gamma_{12}^{2\,(t)}
-\Gamma_{12}^{2\,(t)}\,\Gamma_{12}^{2\,(t)}
+\Gamma_{11}^{2\,(t)}\,\Gamma_{22}^{2\,(t)}\,;\\
K^{(t)}\,g_{12}^{(t)} &=&
\partial_2 \Gamma_{12}^{2\,(t)}
- \partial_1 \Gamma_{22}^{2\,(t)}
-\Gamma_{22}^{1\,(t)}\,\Gamma_{11}^{2\,(t)}
+\Gamma_{12}^{1\,(t)}\,\Gamma_{12}^{2\,(t)}\,,
\rule{0pt}{18pt}
\end{array}
\right.
\end{equation}
with two more equations which can be obtained by interchange of the indices $1$ and $2$. By differentiation w.r.t.\ $t$ at $t=0$ the following four equations result:
\[
\delta K = -\partial_1 X_{12}^{2}\,; \qquad
\partial_2 X_{12}^{2} = 2\,\partial_1 X_{12}^{1}\,;\qquad
\partial_1 X_{12}^{1} = 2\,\partial_2 X_{12}^{2}\,;\qquad
\delta K = -\partial_2 X_{12}^{1}\,.
\]
From the second and the third equation readily follows that both $\partial_2 X_{12}^{2}$ and $\partial_1 X_{12}^{1}$ vanish. From the first and the last equation we see that $\partial_1 (\delta K)$ and $\partial_2 (\delta K)$ vanish. This finishes the proof of the theorem if $K^{(0)}$ vanishes.

\noindent\textbf{(ii.)}---Assume now that the initial curvature $K^{(0)}$ is a strictly positive constant. It is no restriction to assume that $(M,g)$ is the unit sphere, with metric given by $g^{(0)}=\dd u^2+\,(\cos u)^2\,\dd v^2$. The proof of the theorem is more technical in this case because we will need the formulae expressing the variation of the Christoffel symbols in terms of the variation of the metric.

Let us write 
\[
g^{(t)} =
\left\lgroup
\begin{array}{cc}
E^{(t)} & F^{(t)}\\
F^{(t)} & G^{(t)}
\end{array}
\right\rgroup
=
\left\lgroup
\begin{array}{cc}
1 & 0\\
0 & (\cos u)^2
\end{array}
\right\rgroup
+
t
\left\lgroup
\begin{array}{cc}
\delta E\, & \,\delta F\\
\delta F\, & \,\delta G
\end{array}
\right\rgroup
+
\Order(t^2)\,.
\]
A straightforward calculation gives the following expression for the initial Christoffel symbols and their variation:
\[
\begin{array}{rclp{1cm}rcl}
\Gamma_{11}^{1\,(0)} &=& 0\,;              & & X_{11}^{1} &=& 
\frac{1}{2}\delta E_{u}\,; \\
\rule{0pt}{17pt}
\Gamma_{12}^{1\,(0)} &=& 0\,;              & & X_{12}^{1} &=& 
\frac{1}{2}\delta E_{v}+ \tan u\, \delta F\,; \\
\rule{0pt}{17pt}
\Gamma_{22}^{1\,(0)} &=& \cos u \sin u\,;  & & X_{22}^{1} &=&  
\delta F_v -\sin u \cos u\, \delta E - \frac{1}{2}\delta G_u\,;\\
\rule{0pt}{17pt}
\Gamma_{11}^{2\,(0)} &=& 0\,;              & & X_{11}^{2} &=&  
\frac{1}{(\cos u)^2}\left(\delta F_u-\frac{1}{2}\delta E_v\right)\,;\\
\rule{0pt}{17pt}
\Gamma_{12}^{2\,(0)} &=& -\tan u\,;        & & X_{12}^{2} &=&  
\frac{\sin u}{(\cos u)^3}\delta G + \frac{1}{2(\cos u)^2}\delta G_u\,;\\
\rule{0pt}{17pt}
\Gamma_{22}^{2\,(0)} &=& 0\,;              & & X_{22}^{2} &=&   
-\tan u\, \delta F +\frac{1}{2(\cos u)^2}\delta G_v\,.
\end{array}
\]
Here we have denoted $E_v$ for $\partial_2 E$ and we remark that $\delta E_v = \partial_2 \delta E$.
The fact that $\left\{g^{(t)}\right\}$ is an infinitesimally cogeodesical deformation precisely means that (\ref{eq:XX}) is satisfied. This can be rewritten as
\begin{equation}
\label{eq:tsja}
\left\{
\begin{array}{lp{1cm}l}
\textrm{(i.)} &&
\frac{1}{2}\delta E_u = 2\frac{\sin u}{(\cos u)^3}\delta G + \frac{1}{(\cos u)^2}\delta G_u\,;\\
\rule{0pt}{17pt}
\textrm{(ii.)} &&
2\delta F_u = \delta E_v\,;\\
\rule{0pt}{17pt}
\textrm{(iii.)} &&
\delta F_v - \sin u \cos u \, \delta E - \frac{1}{2}\delta G_u =0\,;\\
\rule{0pt}{17pt}
\textrm{(iv.)} &&
\frac{1}{2(\cos u)^2}\delta G_v = \delta E_v + 3 \tan u\, \delta F\,.
\end{array}
\right.
\end{equation}
By applying the operator $\delta$ to both sides of each of the four equations (\ref{eq:Kgij}), which are generally valid, the following information results:
\begin{equation}
\label{eq:deltaKgij}
\left\{
\begin{array}{lp{1cm}l}
\textrm{(i.)} &&
\delta K + \delta E = - \frac{1}{2}\partial_1 X_{11}^{1}\,;\\
\rule{0pt}{17pt}
\textrm{(ii.)} &&
\delta F = \frac{1}{2} \partial_2 X_{11}^{1}-\partial_1 X_{22}^{2}
-\frac{1}{2} \tan u\, X_{22}^{2} \,;\\
\rule{0pt}{17pt}
\textrm{(iii.)} &&
\delta F = \frac{1}{2} \partial_1 X_{22}^{2}-\partial_2 X_{11}^{1}
-\frac{1}{2} \tan u\, X_{22}^{2} \,;\\
\rule{0pt}{17pt}
\textrm{(iv.)} &&
\textrm{(another not so important equation)}.
\end{array}
\right.
\end{equation}
A combination of  (\ref{eq:deltaKgij}.ii.)+2$\times$(\ref{eq:deltaKgij}.iii.), the expressions for $X_{ij}^{k}$ and (\ref{eq:tsja}) results in the following ordinary differential equation for $\delta F$ w.r.t.\ the variable $u$:
\begin{equation}
\label{eq:odedeltaF}
\delta F = -\sin u \cos u \,\delta F_u - \frac{1}{2} (\cos u)^2 \delta F_{uu}\,.
\end{equation}
From (\ref{eq:deltaKgij}.i.) and (\ref{eq:XX}) follows that 
\begin{flalign*}
\qquad\qquad
\delta K_v
&= -\delta E_v -\frac{1}{2} \partial_2\partial_1 X_{11}^{1}\,,&&
\intertext{which can be rewritten as follows, if the expression for $X_{11}^{1}$ and the equations (\ref{eq:tsja}.ii.) and (\ref{eq:odedeltaF}) are subsequently used:}
&= -\delta E_v -\frac{1}{4} \delta E_{uuv}
= -2\delta F_u -\frac{1}{2} \delta F_{uuu}
=0\,.&&
\end{flalign*}
Since the co-ordinates can always be rotated around any point of $ M$, the above equation implies that the directional derivative of $\delta K$ along \textit{any} vector field vanishes.

\noindent\textbf{(iii.)}---The proof is omitted for the case of constant strictly negative curvature.
\end{proof}

\subsection{The Infinitesimal Beltrami Theorem for Arbitrary Dimensions.}
\label{sec:arb_dim}

Let $\left\{g^{(t)}\right\}$ be an infinitesimal deformation of an $n$-dimensional Riemannian manifold $(M,g)$. Besides the symmetric operator $\sigma$ for which $\delta g = \sigma\lrcorner\,g$ holds and the 
variation tensor $X$ of the Levi-Civita connection, we will adopt the following notation:
\begin{equation}
\label{eq:alphalambda}
\alpha = \sigma - \frac{1}{(n+1)}(\textrm{trace}\,\sigma)\,\textrm{id}
\,.
\end{equation}
A characterisation of infinitesimally cogeodesical deformations by means of the tensor $X$ has already been given in lemma \ref{lem:beltrami}. The next lemma, which can be seen as a mere rephrasing of the previous lemma because 
of the relation (\ref{eq:Xsigma}) between the tensors $X$ and $\sigma$,
gives a similar characterisation by means of the tensor $\alpha$.

\begin{lemma_}
\label{lem:beltrami2}
A family of metrics $\left\{g^{(t)}\right\}$ of an $n$-dimensional Riemannian manifold $(M,g)$ (with $g^{(0)}=g$) is an infinitesimally cogeodesical deformation if and only if the following equation is satisfied for all $V,W,Z\in\mathfrak{X}(M)$:
\begin{equation}
\label{eq:alpha}
g((\nabla_Z \alpha)V,W)=\frac{1}{2}V[\mathrm{trace}\,\alpha]\,g(Z,W)
+\frac{1}{2}W[\mathrm{trace}\,\alpha]\,g(Z,V)\,.
\end{equation}
\end{lemma_}

\begin{theorem_}[``\textit{Infinitesimal Beltrami theorem}'']
\label{thm:inf_beltrami_arbdim}
The variation of the sectional curvature of a two-dimensional tangent plane of a space of constant Riemannian curvature under an infinitesimally cogeodesical deformation is a constant which depends neither on the footpoint of the tangent plane, nor on its direction.
\end{theorem_}
\begin{proof}
Assume that $\left\{g^{(t)}\right\}$ is an infinitesimally cogeodesical deformation of a space of constant Riemannian curvature $(M,g)$ of dimension $n\geqslant 2$. According to lemma \ref{lem:beltrami2}, the tensor $\alpha$ as defined in (\ref{eq:alphalambda}) satisfies the linear equation (\ref{eq:alpha}), and consequently the same applies for the tensor 
\[
L^{(t)}= \textrm{id}-t\,\alpha\,,
\]
and this for every $t$. 
Define a one-parameter family of metrics  by the formula
\[
\widetilde{g}^{(t)} = \frac{1}{(\det\,L^{(t)})} \left( \left(L^{(t)}\right)^{-1}\lrcorner\,g\right).
\]
Riemannian invariants which have been constructed w.r.t.\ this family of metrics will bear an extra tilde (\,$\widetilde{\ }$\,) in their notation.
It can readily be checked that these families agree up to first order:
\begin{equation}
\label{eq:g_gtilde}
g^{(t)} = \widetilde{g}^{(t)} + \Order(t^2)\,.
\end{equation}
As a consequence of the circumstance that the tensor $L^{(t)}$ satisfies the equation (\ref{eq:alpha}), the metrics $\widetilde{g}^{(t)}$ and $g$ share their geodesics. This fact has been mentioned in, \textit{e.g.},
\cite{bolsinovmatveev2003}, thm.\ 2;
\cite{matveev2007}, thm.\ 2;
\cite{kiosakmatveev2009}, \S\,2.2 (in which also reference to some classical sources is given).

Thus, whereas $\left\{g^{(t)}\right\}$ was merely known to be an infinitesimally cogeodesical deformation, the newly defined family $\left\{\widetilde{g}^{(t)}\right\}$ of metrics provides us with a cogeodesical deformation (as defined in definition \ref{def:inf_proj_def}). As a consequence of the classical version of Beltrami's theorem, for every $t$, the metric $\widetilde{g}^{(t)}$ has constant sectional curvature (say, $\widetilde{C}^{(t)}$). Define now the constant $\delta\widetilde{C}=\deltatnought  \widetilde{C}^{(t)}$ and let $\widetilde{K}^{(t)}(\Pi)$ stand for the sectional curvature of $(M, \widetilde{g}^{(t)})$ along an arbitrary two-dimensional tangent plane $\Pi$ of $M$. Of course there holds $\widetilde{K}^{(t)}(\Pi)=\widetilde{C}^{(t)}$ and consequently $\deltatnought\widetilde{K}^{(t)}(\Pi)=\delta\widetilde{C}$.

Because of (\ref{eq:g_gtilde}), we have ${K}^{(t)}(\Pi)=\widetilde{K}^{(t)}(\Pi)+\Order(t^2)$ and hence there also holds
$\deltatnought K^{(t)}(\Pi)=\delta\widetilde{C}$, independent of the choice of $\Pi$. This finishes the proof.
\end{proof}

\subsection{A Bibliographical Comment.}

Apparently, the concept of  ``\textit{infinitesimally cogeodesical deformations}'' which was defined above has not been given previously. However, it should be mentioned that a very similar concept has been defined in \S\,2 of \cite{gavr_etal2004} (see also \cite{gavr_1989}) for infinitesimal deformations of Riemannian submanifolds.  
Unfortunately, I find the definition which appears in that article rather ungeometrical and not very precise, because we should drop terms of order $\varepsilon^2$ in the statement of the definition although there does not occur any $\varepsilon$.
As such it appears that equation (7) of that article \cite{gavr_etal2004}, which is perhaps similar to our equations (\ref{eq:X},\ref{eq:alpha}), is taken as the defining equation and the starting point for the study of such deformations in 
\cite{gavr_etal2004}.

\vspace{5mm}
\addcontentsline{toc}{section}{Acknowledgements}
\noindent\textbf{Acknowledgements.}
I would like to express my gratitude for useful dicussions towards Professors J.\ Mike\v{s} and V.S.\ Matveev, the latter of whom has kindly suggested some crucial elements of the proof of theorem~\ref{thm:inf_beltrami_arbdim}.

The author, who was employed at K.U.Leuven during the commencement of this work, while he was supported by Masaryk University (Brno) during its conclusion, is thankful to both these institutions. 
This research was partially supported by 
the Research Foundation Flanders (project G.0432.07) and the Eduard \v{C}ech Center for Algebra and Geometry (Basic Research Center no. LC505).

\bibliographystyle{amsplain}


\begin{center}
\rule{0.6\textwidth}{1.0pt}
\end{center}

\end{document}